\documentclass[12pt,reqno]{amsart}
\usepackage{amssymb}
\usepackage{amsfonts}
\usepackage{graphicx}
\usepackage{color}
\usepackage{hyperref}

\def\N{\mathbb{N}}

\begin{document}

\title[Matrix products of binomial coefficients$\ldots$]{Matrix products of binomial coefficients and unsigned Stirling numbers}

\author[M.~Kne\v{z}evi\'{c}]{Marin Kne\v{z}evi\'{c}}

\email{marin.k59@gmail.com}

\author[V.~Kr\v{c}adinac]{Vedran Kr\v{c}adinac}

\thanks{The second author has been supported by the Croatian Science Foundation
under the project $6732$.}

\email{vedran.krcadinac@math.hr}

\address{Department of Mathematics, Faculty of Science, University of Zagreb,
Bijeni\v{c}ka~30, HR-10000 Zagreb, Croatia}

\author[L.~Reli\'{c}]{Lucija Reli\'{c}}

\email{lucijarelic7@gmail.com}

\dedicatory{Dedicated to Professor Dragutin Svrtan on the occasion
of his 70$^{th}$ birthday.}


\begin{abstract}
We study sums of the form $\sum_{k=m}^n a_{nk} b_{km}$, where
$a_{nk}$ and $b_{km}$ are binomial coefficients or unsigned Stirling
numbers. In a few cases they can be written in closed form. Failing
that, the sums still share many common features: combinatorial
interpretations, Pascal-like recurrences, inverse relations with
their signed versions, and interpretations as coefficients of change
between polynomial bases.
\end{abstract}

\date{December 30, 2020}

\maketitle

\section{Introduction}

Let $a_{nk}$ and $b_{nk}$ be binomial coefficients ${n\choose k}$ or
unsigned Stirling numbers of the first and second kind:
$\left[n\atop k\right]$ and $\left\{n\atop k\right\}$ in
Karamata-Knuth notation, respectively. We are interested in the sum
$\sum_{k=m}^n a_{nk} b_{km}$. Denote by $A=[a_{nk}]$ and
$B=[b_{nk}]$ the corresponding infinite lower-triangular matrices,
indexed by the non-negative integers. The sum can be interpreted as
the $(n,m)$-entry of the matrix product $A\cdot B$.

Our motivation are the following two sums, that can be written in
closed form:
\begin{equation}\label{cfbins2}
\sum_{k=m}^n {n\choose k} \left\{k\atop m\right\} = \left\{n+1\atop
m+1\right\},
\end{equation}
\begin{equation}\label{cfs1bin}
\sum_{k=m}^n \left[\rule{0mm}{3mm}n\atop k\right]{k\choose m} =
\left[n+1\atop m+1\right].
\end{equation}
These are identities $(6.15)$ and $(6.16)$ in the book
\emph{Concrete mathematics}~\cite{GKP94}. Furthermore, the following
sum are the Lah numbers~\cite{IL54, IL55}, denoted by $L(n,m)$:
\begin{equation}\label{cflah}
\sum_{k=m}^n \left[\rule{0mm}{3mm}n\atop k\right]\left\{k\atop
m\right\}  = \frac{n!}{m!}{n-1\choose m-1}.
\end{equation}
Closed-form expressions for two more sums are given in the next
section. In the remaining cases we do not know closed forms, but the
sums share many common features. Combinatorial interpretations are
outlined in Section~2. Row sums of the matrix $A\cdot B$ are
discussed in Section~3. Recurrences similar to Pascal's formula are
proved in Section~4. Inverse relations with signed versions of the
sums are given in Section~4, and in Section~5 the sums are
interpreted as coefficients of change between various polynomial
bases.

\begin{table}
\begin{center}
\begin{tabular}{|cc|cc|cc|}
\hline
\rule{0mm}{6mm}$\sum_k {n\choose k}{k\choose m}$ & A038207 &
$\sum_k {n\choose k}\left[k\atop m\right]$ & A094816 &
$\sum_k {n\choose k}\left\{k\atop m\right\}$ & A008277\\[2mm]
\hline
\rule{0mm}{6mm}$\sum_k \left[n\atop k\right]{k\choose m}$ & A130534 &
$\sum_k \left[n\atop k\right]\left[k\atop m\right]$ & A325872$^*$ &
$\sum_k \left[n\atop k\right]\left\{k\atop m\right\}$ & A271703 \\[2mm]
\hline
\rule{0mm}{6mm}$\sum_k \left\{n\atop k\right\}{k\choose m}$ & A049020 &
$\sum_k \left\{n\atop k\right\}\left[k\atop m\right]$ & A129062 &
$\sum_k \left\{n\atop k\right\}\left\{k\atop m\right\}$ & A130191 \\[2mm]
\hline
\end{tabular}
\vskip 4mm
\end{center}
\caption{References for the sums in the On-line encyclopedia of
integer sequences~\cite{oeis}.}\label{table1}
\end{table}

Our main focus are sums with $a_{nk}$ and $b_{nk}$ equal to
${n\choose k}$, $\left[n\atop k\right]$ or $\left\{n \atop
k\right\}$, but we include some results on sums with $L(n,k)$ as
well. Our nine main sums are listed in the On-line encyclopedia of
integer sequences~\cite{oeis} as ``triangles read by rows'', see
Table~\ref{table1} (the entry marked by $^*$ is a signed version).
Many properties of the sums are reported in~\cite{oeis}. We repeat
some properties from~\cite{oeis}, and reveal some new properties of
the sums. We mainly deal with properties that can be proved by
counting arguments.

\section{Combinatorial interpretations and closed forms}

Identity~\eqref{cfbins2} follows from the usual combinatorial
interpretation of Stirling numbers of the second kind. The
right-hand side $\left\{n+1\atop m+1\right\}$ counts partitions of
an $(n+1)$-element set~$S$ into $m+1$ blocks. The left-hand side is
obtained by distinguishing an element of~$S$ and assuming it is
covered by a block containing $n-k$ other elements of~$S$.

The Stirling number of the first kind~$\left[n\atop m\right]$ is
usually interpreted as the number of permutations of degree~$n$ with
exactly~$m$ cycles. For the proof of identity~\eqref{cfs1bin},
however, another combinatorial interpretation is more useful. Let
$\Delta_n$ denote the board that remains after removing all tiles
below and on the side diagonal of a $n \times n$ square board.
Consider a new chess piece similar to a rook, but that only attacks
tiles in its row. Loehr~\cite{L11} calls this piece a weak rook or
wrook. He goes on to show that the number of ways of placing $k$
identical wrooks on $\Delta_n$ so that they don't attack each other
is equal to $\left[n\atop n-k\right]$. The right-hand side
of~\eqref{cfs1bin} can now be interpreted as the number of ways of
placing $n-m$ non-attacking wrooks on $\Delta_{n+1}$. The left-hand
side is obtained by partitioning the set of placements based on the
number of wrooks not in the first column, but on the remaining
$\Delta_n$ subboard. If $n-k$ wrooks are placed on the $\Delta_n$
subboard, the remaining $k-m$ wrooks must be placed in $k$ rows of
the first column that are not attacked. This can be done in
$\left[n\atop k\right]{k\choose k-m} = \left[n\atop
k\right]{k\choose m}$ ways.

In~\cite{PP07}, a nice combinatorial interpretation of the Lah
numbers is given: $L(n,m)$ counts partitions of an $n$-element set
into $m$ lists, i.e.\ non-empty totally ordered subsets. From this,
identity~\eqref{cflah} follows easily. Two more closed-form
identities are
\begin{equation}\label{cfbinbin}
\sum_{k=m}^n {n\choose k}{k\choose m}  = 2^{n-m} {n\choose m}
\end{equation}
and
\begin{equation}\label{cflahlah}
\sum_{k=m}^n L(n,k)\,L(k,m)  = 2^{n-m} L(n,m).
\end{equation}
The left-hand side of~\eqref{cfbinbin} counts pairs $(K,M)$, where
$K\supseteq M$ are subsets of a fixed $n$-element set~$S$, and
$|M|=m$. The right-hand side is obtained by choosing~$M$ first.
Identity~\eqref{cflahlah} has an analogous proof using the
combinatorial interpretation of Lah numbers from~\cite{PP07}. For
two partitions $\Pi_1$, $\Pi_2$ of~$S$ into lists, we write
$\Pi_1\le \Pi_2$ and say that $\Pi_1$ is a refinement of $\Pi_2$ if
each list in~$\Pi_1$ is a sublist of some list in~$\Pi_2$. The
left-hand side of~\eqref{cflahlah} counts pairs $(\Pi_1,\Pi_2)$ with
$\Pi_1\le \Pi_2$ and $|\Pi_2|=m$. We first partition~$S$
into~$\Pi_1$ with $|\Pi_1|=k\ge m$, then partition these $k$ lists
into~$m$ lists containing lists of~$\Pi_1$ as elements, and finally
concatenate the list of lists into elements of~$\Pi_2$. We get the
right-hand side by partitioning~$S$ into~$m$ lists of~$\Pi_2$ first
and splitting them up into smaller lists of~$\Pi_1$. This can be
done in $2^{n-m}$ ways: breaks can be made in $n-m$ places, before
every element except at the beginnings of the lists of~$\Pi_2$.

Similar combinatorial interpretations can be given for the double
Stirling sums $\sum_{k=m}^n \left\{n\atop k\right\} \left\{k\atop
m\right\}$ and $\sum_{k=m}^n \left[ n\atop k\right] \left[ k\atop
m\right ]$. The former is the number of pairs $(\Pi_1,\Pi_2)$ of
ordinary partitions of an $n$-element set into blocks (non-empty
subsets), with $\Pi_1$ a refinement of $\Pi_2$ and $|\Pi_2|=m$. The
latter is the number of pairs of permutations $(\pi_1,\pi_2)$, where
$\pi_1$ is of degree~$n$, $\pi_2$ is of degree equal to the number
of cycles of~$\pi_1$, and~$\pi_2$ has exactly~$m$ cycles. However,
since we do not see an easy way to count the pairs if $\Pi_2$ or
$\pi_2$ are chosen first, this does not lead to closed-form
expressions for these two sums.

Combinatorial interpretations of the remaining sums in
Table~\ref{table1} are as follows. The sum $\sum_{k=m}^n
\left\{n\atop k\right\} {k\choose m}$ is the number of ordinary
partitions of an $n$-element set with $m$ blocks colored red, and
the other blocks colored blue. The sum $\sum_{k=m}^n {n\choose k}
\left[ k\atop m\right]$ counts permutations with exactly~$m$ cycles
of all subsets of an $n$-element set. Finally, $\sum_{k=m}^n
\left\{n\atop k\right\} \left[k\atop m\right]$ is the number of
pairs $(\Pi,\pi)$, where~$\Pi$ is a partition of an $n$-element set
into blocks, and~$\pi$ is a permutation with exactly~$m$ cycles of
degree equal to the number of blocks of~$\Pi$. Again, this does not
lead to closed-form expressions, but the combinatorial
interpretations of the sums will be used in the following sections
to prove their properties.

\section{Row sums}

From the combinatorial interpretation of binomial coefficients and
Stirling numbers, it is clear that $\sum_{m=0}^n {n\choose m}= 2^n$,
$\sum_{m=0}^n \left[ n\atop m \right]=n!$, and $\sum_{m=0}^n \left\{
n\atop m \right\} = B_n$. Here $B_n$ is the $n$-th Bell number,
i.e.\ the total number of partitions of an $n$-element set into
blocks. In this section the goal is to determine
\[ \sum_{m=0}^n \sum_{k=m}^n a_{nk} b_{km}. \]
This is the sum of the $n$-th row of the matrix $A\cdot B$.

From identities~\eqref{cfbins2} and~\eqref{cfs1bin}, we have
\[ \sum_{m=0}^n \sum_{k=m}^n {n\choose k} \left\{k\atop m\right\}
= B_{n+1} \,\,\,\mbox{ and }\,\,\, \sum_{m=0}^n \sum_{k=m}^n
\left[\rule{0mm}{3mm}n\atop k\right]{k\choose m} = (n+1)!. \]
Another row sum that can be written in closed form is
\[ \sum_{m=0}^n \sum_{k=m}^n {n\choose k}{k\choose m} = 3^n. \]
The left-hand side counts the total number of pairs $(K,M)$ of
subsets of an $n$-set $S$ with $K\supseteq M$. An alternative way of
counting is to decide for each element $x\in S$ whether it is in
$M$, in $K\setminus M$, or in $S\setminus K$, leading to the
right-hand side.

The row sum of Lah numbers
\[ \sum_{m=0}^n \sum_{k=m}^n \left[\rule{0mm}{3mm}n\atop
k\right]\left\{k\atop m\right\} = \sum_{m=0}^n L(n,m) \] can be
interpreted as the total number of partitions of an $n$-set into
lists. This is the ``Lah equivalent'' of the Bell number~$B_n$ and
we will denote it by~$L_n$. In~\cite{oeis}, this is sequence number
A000262.

Row sums of the remaining matrix products are given in the sequel.
They have nice combinatorial interpretations and can be simplified
to single sums.
\[ \sum_{m=0}^n \sum_{k=m}^n \left\{\rule{0mm}{3mm}n\atop k\right\}
\left[k\atop m\right] = \sum_{m=0}^n \left\{\rule{0mm}{3mm}n\atop
m\right\} m!.\] These are the Fubini numbers, sequence A000670
in~\cite{oeis}. They count ordered partitions of an $n$-set, or weak
orders on $n$ elements.\\
\[ \sum_{m=0}^n \sum_{k=m}^n \left[\rule{0mm}{3mm}n\atop k\right]
\left[k\atop m\right] = \sum_{m=0}^n \left[\rule{0mm}{3mm}n\atop
m\right] m!.\] This is the number of ordered factorizations of
permutations of degree~$n$ into cycles, sequence A007840
in~\cite{oeis}.\\
\[ \sum_{m=0}^n \sum_{k=m}^n \left\{\rule{0mm}{3mm}n\atop k\right\}
{k\choose m} = \sum_{m=0}^n \left\{\rule{0mm}{3mm}n\atop m\right\}
2^m. \] The total number of partitions of an $n$-set with blocks
colored red or blue. This is sequence A001861 in~\cite{oeis}.\\
\[ \sum_{m=0}^n \sum_{k=m}^n {n\choose k}
\left[k\atop m\right] = \sum_{m=0}^n \frac{n!}{m!}.\] The total
number of lists with elements from an $n$-set. Sequence number
A000522 in~\cite{oeis}.\\
\[ \sum_{m=0}^n \sum_{k=m}^n \left\{\rule{0mm}{3mm}n\atop k\right\}
\left\{k\atop m\right\} = \sum_{m=0}^n \left\{\rule{0mm}{3mm}n\atop
m\right\} B_m. \] The total number of pairs $(\Pi_1,\Pi_2)$ of
partitions of an $n$-set with $\Pi_1\le \Pi_2$. In~\cite{oeis}, this
is sequence number A000258.

\section{Pascal-like recurrences}

The binomial coefficients can be computed by Pascal's formula
\[ {n\choose m}={n-1\choose m-1}+{n-1\choose m},\kern 3mm {n\choose
0}={n\choose n}=1. \] Analogous recurrences for Stirling numbers are
\[ \left[\rule{0mm}{3mm}n\atop m\right]=\left[n-1\atop
m-1\right]+(n-1)\left[n-1\atop m\right], \kern 3mm
\left[\rule{0mm}{3mm}n\atop 0\right]=\delta_{n0},\kern 3mm
\left[\rule{0mm}{3mm}n\atop n\right]=1,\]
\[ \left\{n\atop m\right\}=\left\{n-1\atop m-1\right\}+m\left\{n-1\atop
m\right\},\kern 3mm \left\{\rule{0mm}{3mm}n\atop
0\right\}=\delta_{n0},\kern 3mm \left\{\rule{0mm}{3mm}n\atop
n\right\}=1,\]
and for Lah numbers
\[ L(n,m)=L(n-1,m-1)+(n+m-1)L(n-1,m) \]
with boundary values $L(n,0)=\delta_{n0}$ (the Kronecker delta),
$L(n,n)=1$. See~\cite{PP07} for proofs of the formulae by
distinguishing an element of the underlying $n$-set~$S$ and
counting. Our sums satisfy similar recurrences that can also be
established by counting arguments.

To make the formulae more readable, we denote the sum $\sum_{k=m}^n
a_{nk}\, b_{km}$ by $\left| n\atop m\right|$. For example, the
double binomial sum $\left| n\atop m\right|=\sum_{k=m}^n {n\choose
k}{k\choose m}$ can be computed from
\[ \left|\rule{0mm}{3mm}n\atop
m\right|=\left|n-1\atop m-1\right|+2\left|n-1\atop m\right|,\kern
3mm \left|\rule{0mm}{3mm}n\atop 0\right|=2^n,\kern 3mm
\left|\rule{0mm}{3mm}n\atop n\right|=1. \] This sum also satisfies
the absorption identity $\left|n\atop m\right|= \frac{n}{m}
\left|n-1\atop m-1\right|$, just like the binomial coefficients.

By~\eqref{cfbins2}, the sum $\left|n\atop m\right|=\sum_{k=m}^n
{n\choose k} \left\{k\atop m\right\}$ are shifted Stirling numbers
of the second kind. Therefore,
\[ \left|\rule{0mm}{3mm}n\atop
m\right|=\left|n-1\atop m-1\right|+(m+1)\left|n-1\atop
m\right|,\kern 3mm \left|\rule{0mm}{3mm}n\atop
0\right|=\left|\rule{0mm}{3mm}n\atop n\right|=1. \]

Similarly, by~\eqref{cfs1bin}, the sum $\left|n\atop
m\right|=\sum_{k=m}^n \left[n\atop k\right]{k\choose m}$ satisfies
\[ \left|\rule{0mm}{3mm}n\atop
m\right|=\left|n-1\atop m-1\right|+n\left|n-1\atop m\right|,\kern
3mm \left|\rule{0mm}{3mm}n\atop 0\right|=n!,\kern 3mm
\left|\rule{0mm}{3mm}n\atop n\right|=1. \]

The sum $\left|n\atop m\right|=\sum_{k=m}^n {n\choose k} L(k,m)$ is
sequence A271705 in~\cite{oeis}, where the following recurrence is
given:
\[ \left|\rule{0mm}{3mm}n\atop
m\right|=\frac{n}{m} \left|n-1\atop m-1\right|+n \left|n-1\atop
m\right|,\kern 3mm \left|\rule{0mm}{3mm}n\atop
0\right|=\left|\rule{0mm}{3mm}n\atop n\right|=1. \] The sum
$\left|n\atop m\right|=\sum_{k=m}^n L(n,k){k\choose m}$ is sequence
A059110. It satisfies the same recurrence, but with different
boundary values $\left|n\atop 0\right|=L_n$.

The double Lah sum $\left|n\atop m\right|=\sum_{k=m}^n L(n,k)
L(k,m)$ satisfies
\[ \left|\rule{0mm}{3mm}n\atop
m\right|=\frac{n}{m} \left|n-1\atop m-1\right|+2n \left|n-1\atop
m\right|,\kern 3mm \left|\rule{0mm}{3mm}n\atop
0\right|=\delta_{n0},\kern 3mm \left|\rule{0mm}{3mm}n\atop
n\right|=1. \] This sum also satisfies the absorption-like identity
$\left|n\atop m\right|= \frac{n-m+1}{2m(m-1)} \left|n\atop
m-1\right|$, while the Lah numbers satisfy $L(n,m)=
\frac{n-m+1}{m(m-1)} L(n,m-1)$.

Now let $\left|n\atop m\right|= \sum_{k=m}^n {n\choose k}
\left[k\atop m\right]$. This sum cannot be computed from the two
values $\left|n-1\atop m-1\right|$ and $\left|n-1\atop m\right|$
alone, but we can give a Pascal-like recurrence involving more
previous values:
\[ \left|\rule{0mm}{3mm}n\atop
m\right|=\sum_{k=0}^{n-m} (n-1)^{\underline{k}} \left|n-1-k\atop
m-1\right|+ \left|n-1\atop m\right|,\kern 3mm
\left|\rule{0mm}{3mm}n\atop 0\right|= \left|\rule{0mm}{3mm}n\atop
n\right|=1. \] Here $(n-1)^{\underline{k}} = (n-1)\cdot(n-2)\cdots
(n-k)$ is the falling factorial. For the proof, recall that
$\left|n\atop m\right|$ counts permutations with $m$ cycles of
subsets $T\subseteq S$, where $S$ is a set of $n$ elements. Fix an
element $x\in S$ and divide the permutations according to whether
they contain~$x$ or do not contain~$x$. In the latter case there are
clearly $\left|n-1\atop m\right|$ permutations. In the former case,
assume $x$ is in a cycle with $k$ other elements of~$S$. This cycle
can be chosen in $(n-1)^{\underline{k}}$ ways, and $m-1$ more cycles
have to be chosen from the remaining $n-1-k$ elements. Thus, there
are $\sum_{k=0}^{n-m} (n-1)^{\underline{k}} \left|n-1-k\atop
m-1\right|$ permutations containing~$x$.

The sum $\left|n\atop m\right|= \sum_{k=m}^n \left\{n\atop k\right\}
{k\choose m}$ counts partitions of $S$ with $m$ blocks colored red,
and the other blocks colored blue. Again, fix an element $x\in S$.
If $x$ is in a red block alone, there are $\left|n-1\atop
m-1\right|$ partitions. If $x$ is in a red block with some other
elements of~$S$, there are $m \left|n-1\atop m\right|$ partitions.
Finally, if $x$ is in a blue block with $k$ other elements of $S$,
there are ${n-1\choose k}\left|n-1-k\atop m\right|$ such partitions.
Therefore, the following recursion holds:
\[ \left|\rule{0mm}{3mm}n\atop
m\right|= \left|n-1\atop m-1\right| + m \left|n-1\atop m\right| +
\sum_{k=0}^{n-m-1} {n-1\choose k}\left|n-1-k\atop m\right|. \] The
boundary values are $\left|n\atop 0\right|=B_n$ and $\left|n\atop
n\right|=1$.

The sum $\left|n\atop m\right|= \sum_{k=m}^n \left\{n\atop k\right\}
\left\{k\atop m\right\}$ counts pairs $(\Pi_1,\Pi_2)$ of partitions
of~$S$, where $\Pi_2$ has $m$ blocks and $\Pi_1$ is a refinement of
$\Pi_2$. Now let the fixed element $x\in S$ be contained in a block
of~$\Pi_2$ of size~$k$. We can choose this block in ${n-1\choose
k-1}$ ways and partition it into blocks of $\Pi_1$ in $B_k$ ways.
The remaining blocks of $\Pi_2$ and $\Pi_1$ can be chosen in
$\left|n-k\atop m-1\right|$ ways. Therefore,
\[ \left|n\atop m\right|=\sum\limits_{k=1}^{n-m+1} {n-1\choose k-1} B_k
\left|n-k\atop m-1\right|,\kern 3mm \left|\rule{0mm}{3mm}n\atop
0\right|=\delta_{n0},\kern 3mm \left|\rule{0mm}{3mm}n\atop
n\right|=1.  \]

For the sum $\left|n\atop m\right|= \sum_{k=m}^n \left[n\atop
k\right] \left[k\atop m\right]$ a similar argument leads to the
recurrence
\[ \left|n\atop m\right|=\sum\limits_{k=1}^{n-m+1} {n-1\choose k-1}
\sum_{i=1}^k \left[k\atop i\right] (i-1)! \left|n-k\atop
m-1\right|,\kern 3mm \left|\rule{0mm}{3mm}n\atop
0\right|=\delta_{n0},\kern 3mm \left|\rule{0mm}{3mm}n\atop
n\right|=1,  \] and for the sum $\left|n\atop m\right|= \sum_{k=m}^n
\left\{n\atop k\right\} \left[k\atop m\right]$ to the recurrence
\[ \left|n\atop m\right|=\sum\limits_{k=1}^{n-m+1} {n-1\choose k-1}
\sum_{i=1}^k \left\{k\atop i\right\} (i-1)! \left|n-k\atop
m-1\right|,\kern 3mm \left|\rule{0mm}{3mm}n\atop
0\right|=\delta_{n0},\kern 3mm \left|\rule{0mm}{3mm}n\atop
n\right|=1.  \] However, these increasingly complex recurrences
become less useful as the coefficients are more difficult to
evaluate than the sum $\left|n\atop m\right|$ directly.

\section{Inverses}

We denote signed versions of the Stirling numbers and their
relatives by an exponent $^\sigma$, e.g.\ $\left[n\atop
m\right]^\sigma = (-1)^{n-m} \left[n\atop m\right]$. For the matrix
$A=[a_{nm}]$, we denote $A^\sigma=[a_{nm}^\sigma]=
[(-1)^{n-m}a_{nm}]$. To avoid excessive bracketing, we write
$\left[n\atop m\right]^{-1}$ for the inverse matrix $A^{-1}$.

It is well-known that ${n\choose m}^{-1}={n\choose m}^\sigma$,
$\left[n\atop m\right]^{-1}=\left\{n\atop m\right\}^\sigma$, and
$\left\{n\atop m\right\}^{-1}=\left[n\atop m\right]^\sigma$. From
this and the properties $(A\cdot B)^{-1}=B^{-1}\cdot A^{-1}$ and
$(A\cdot B)^\sigma=A^\sigma\cdot B^\sigma$, we can determine
inverses of our sums. For example, let $a_{mn}=\left[n\atop
m\right]$ and $b_{mn}=\left\{n\atop m\right\}$. Then $A\cdot B$ are
the Lah numbers $L(n,m)$ and we have
\[ L(n,m)^{-1}=(A\cdot
B)^{-1}=B^{-1}\cdot A^{-1}=A^\sigma\cdot B^\sigma=(A\cdot B)^\sigma
= L(n,m)^\sigma. \] Similarly it follows
\[ \left(\sum_{k=m}^n {n\choose k} \left\{k\atop
m\right\}\right)^{-1} = \left(\sum_{k=m}^n
\left[\rule{0mm}{3mm}n\atop k\right]{k\choose m}\right)^\sigma, \]
\[ \left(\sum_{k=m}^n \left\{\rule{0mm}{3mm}n\atop
k\right\} {k\choose m}\right)^{-1} = \left(\sum_{k=m}^n {n\choose k}
\left[k\atop m\right]\right)^\sigma, \]
\[ \left(\sum_{k=m}^n \left\{\rule{0mm}{3mm}n\atop
k\right\}\left\{k\atop m\right\}\right)^{-1} = \left(\sum_{k=m}^n
\left[\rule{0mm}{3mm}n\atop k\right]\left[k\atop
m\right]\right)^\sigma, \] and
\[ \left(\sum_{k=m}^n \left\{\rule{0mm}{3mm}n\atop
k\right\}\left[k\atop m\right]\right)^{-1} = \left(\sum_{k=m}^n
\left\{\rule{0mm}{3mm}n\atop k\right\}\left[k\atop
m\right]\right)^\sigma. \] The exponents $^{-1}$ and $^\sigma$ can
be exchanged in the formulae above.

\section{Polynomial bases}

We denote the falling factorials by $x^{\underline{n}}=x(x-1)\cdots
(x-n+1)$ and the rising factorials by $x^{\overline{n}}=x(x+1)\cdots
(x+n-1)$, following~\cite{GKP94}. It is known that Stirling numbers
of the second kind are coefficients of change from the standard
polynomial basis of powers $(x^n \mid n\in \N_0)$ to the basis of
falling factorials $(x^{\underline{m}} \mid m\in\N_0)$:
\begin{equation}\label{powtoff}
x^n = \sum_{m=0}^n \left\{n \atop m\right\} x^{\underline{m}}.
\end{equation}
Stirling numbers of the first kind are coefficients of change from
$(x^{\overline{n}})$ to $(x^m)$:
\begin{equation}\label{rftopow}
x^{\overline{n}} = \sum_{m=0}^n \left[n \atop m\right] x^m.
\end{equation}
From~\eqref{powtoff} and~\eqref{rftopow} it follows that
coefficients of change from $(x^{\overline{n}})$ to
$(x^{\underline{m}})$ are the Lah numbers $L(n,m)=\sum_{k=m}^n
\left[n \atop k\right] \left\{k \atop m\right\}$:
\[ x^{\overline{n}} = \sum_{m=0}^n L(n,m) x^{\underline{m}}.
\]
Ivo Lah's original definition of his numbers~\cite{IL54, IL55} was a
signed version of this relation. We concentrate on changes between
polynomial bases with non-negative coefficients. The opposite
changes have inverse coefficients, with alternating signs as shown
in the previous section.

From the binomial theorem $(1+x)^n = \sum_{m=0}^n {n\choose m} x^m$
and~\eqref{powtoff} we see that the sums $\sum_{k=m}^n {n\choose k}
\left\{k\atop m\right\}$ can be interpreted as coefficients of
change from the polynomial basis $((1+x)^n)$ to the basis of falling
factorials $(x^{\underline{m}})$:
\[ (1+x)^n = \sum_{m=0}^n \left( \sum_{k=m}^n {n\choose k}
\left\{k\atop m\right\} \right) x^{\underline{m}}. \] Similar
interpretations can be given to other sums from Table~\ref{table1}.
The double binomial sums $\sum_{k=m}^n {n\choose k}{k\choose m}$ are
coefficients of change from the basis $((2+x)^n)$ to the standard
basis $(x^m)$. The sums $\sum_{k=m}^n {n\choose k}\left[k\atop
m\right]$ are coefficients of change from the basis $(\sum_{k=0}^n
{n\choose k} x^{\overline{k}})$ to the standard basis. The former
basis contains a special case of Charlier polynomials~\cite{oeis}, a
family of orthogonal polynomials that can be written in terms of the
generalized hypergeometric function (see~\cite{OED15}). Double
Stirling sums of the second kind $\sum_{k=m}^n \left\{n\atop
k\right\}\left\{k\atop m\right\}$ are coefficients of change from
the basis of Bell polynomials $(B_n(x))$, $B_n(x)=\sum_{k=0}^n
\left\{n\atop k\right\} x^k$ to the basis of falling factorials
$(x^{\underline{m}})$.

The two families of sums $\sum_{k=m}^n \left\{n\atop
k\right\}\left[k\atop m\right]$ and $\sum_{k=m}^n \left[n\atop
k\right]\left[k\atop m\right]$ can be seen as coefficients of the
polynomials $\sum_{k=0}^n \left\{n\atop k\right\} x^{\overline{k}}$
and $\sum_{k=0}^n \left[n\atop k\right] x^{\overline{k}}$, i.e.\
coefficients of change to the standard basis $(x^m)$. Similarly, the
two families $\sum_{k=m}^n \left\{n\atop k\right\}{k\choose m}$ and
$\sum_{k=m}^n \left[n\atop k\right]{k\choose m}$ are coefficients of
Bell polynomials $B_n(1+x)$ and polynomials $\sum_{k=0}^n
\left[n\atop k\right] (1+x)^k$, respectively.

\end{document}